\documentclass{article}
\usepackage{amssymb}

\usepackage{graphicx}
\usepackage{amsmath}

\newtheorem{theorem}{Theorem}

\begin{document}

\title{Bilaplacian on a Riemannian manifold and Levi-Civita connection}
\author{R\'emi L\'eandre   \\
Laboratoire  de Math\'ematiques, Universit\'e de Franche-Comt\'e,
\\ 25030, Besan\c  con, FRANCE.
\\email: remi.leandre@univ-fcomte.fr
} \maketitle

\begin{abstract} We generalize for the Bilaplacian the Eells-Elworthy-Malliavin construction of the Brownian motion on a Riemannian manifold.
\end{abstract} 
\section{Introduction}
It is classical in stochastic analysis that the horizontal lift of a diffusion is useful in order to construct canonically the Brownian motion on a Riemannian manifold ([1], [2]). We extend in this note this classical relation on the Brownian motion and the (degenerated) process associated to the horizontal Laplacian to the case of a Bilaplacian. See [7] in the subelliptic case.

We consider a compact Riemannian oriented manifold $M$ endowed with its normalized Riemannian measure $dx$.  $x$ is the generic element of $M$ which is of dimension $m$.
We consider the special orthonormal frame bundle $SO(M)$ endowed with the Levi-Civita connection with canonical projection $\pi$ on $M$. $u$ is the generic element of $SO(M)$. We consider the canonical vector fields $X_i$ on $SO(M)$ and the associated horizontal Laplacian 
\begin{equation} L = \sum_{i=1}^m X_i^2\end{equation}
$\Delta$ is the Laplace-Beltrami operator on $M$ and $\Delta^2$ the associated Bilaplacian.  The Bilaplacian is elliptic, symmetric on $L^2(dx)$ and by elliptic theory ([5], [6])generates a unique contraction semi group $P_t^{\Delta}$ on $L^2(dx)$ with generic element $f$.

We can glue on the fiber the normalized "Haar" measures (We omitt the details) such that we get a probability measure $du$ on $SO(M)$. $L^2$ is symmetric, densely defined on $L^2(du)$ and therefore admits a self-adjoint extension which generates a contraction semi-group $P_t^L$ on $L^2(du)$.
\begin{theorem}If $f$ is a smooth function on $M$, we have if $\pi u = x$
\begin{equation}P_t^{\Delta}[f](x) = P_t^L[f\circ \pi](u)\end{equation}
\end{theorem}
This theorem enters in our general program to extend stochastic analysis tools to the  general theory of linear semi-group (See [8] and [9] for reviews).
\section{Proof of the theorem}
If $u_0$ belongs to $SO(\mathbb{R}^m)$ we get clearly
\begin{equation}P_t^L[f\circ \pi](uu_0) = P_t^L[f\circ \pi(.u_0)](u ) = P_t^L[f\circ \pi](u)\end{equation}
Therefore, $u\rightarrow  P_t^L[f\circ \pi](u)$ defines a function on $M$. This function when $t$ is moving defines a semi-group. Namely
\begin{equation}P_{t+s}^L[f\circ \pi](u) = P_t^L[P_s^L[f\circ \pi]](u) =  P_t^L[P_s^L[f\circ \pi]\circ \pi](u) \end{equation}
where is the right-hand side $ P_s^L[f\circ \pi]$ is seen as a function on $M$.

We are therefore in presence of two semi-groups on $L^2(dx)$. But
\begin{equation}L(f\circ \pi) = (\Delta f) \circ \pi\end{equation}
such that 
\begin{equation}L^2(f\circ \pi) = (\Delta^2f)\circ \pi\end{equation}
Since there is only one semi-group generated by $\Delta^2$, the result holds.

$\diamondsuit$

\end{document}